\documentclass[11pt]{amsart}
\usepackage{amsmath,amssymb,amscd}
\usepackage{amsmath,amssymb}
\usepackage{mathrsfs}
\usepackage{color}
\usepackage{amsxtra, amsmath}
\usepackage{amssymb, amscd}
\usepackage{graphicx}
\usepackage{mathtools}
\usepackage{cite}

\theoremstyle{plain}
\newtheorem{thm}{Theorem}[section]

\theoremstyle{definition}

\theoremstyle{remark}
\newtheorem{remark}[thm]{Remark}

\title[Time-and-band limiting for symmetric classical Dunkl orthogonal polynomials]
{Time-and-band limiting for symmetric classical Dunkl orthogonal polynomials}

\author{M. M. Castro \and  F. A.  Gr\"unbaum}

\address{ Departamento de Matem\'atica Aplicada II and IMUS, Escuela Politécnica Superior, Universidad de Sevilla, 41011, Sevilla, Spain}
\email{mirta@us.es}
\thanks{The research of M. M. Castro was partially supported by  PID2024-155593NB-C21 (FEDER(EU)/Ministerio de Ciencia, Innovación y Universidades-Agencia Estatal de Investigación), FQM-262 (Junta de Andalucía) and  IMUS-Mar\'{\i}a de Maeztu grant CEX2024-001517-M (funded by MICIU/AEI/ 10.13039/501100011033). }

\address{Department of Mathematics, University of California, Berkeley
CA 94705, USA}
\email{grunbaum@math.berkeley.edu}



\begin{document}

\begin{abstract}

The problem of ``time-and-band limiting'' has been previously considered for
expansions in orthonormal polynomials that are common eigenfunctions of some second order differential operator. In this paper we consider the situation of orthonormal polynomials that are eigenfunctions of some second order differential-difference operator as introduced by Ch. Dunkl. These polynomials appear already in the work of G. Szeg\H{o} and at the same time play a role in supersymmetric quantum mechanics.

\end{abstract}

\maketitle

\section{Introduction}\label{intro}

	The celebrated series of papers by D. Slepian, H. Landau and H. Pollak, at Bell Labs in the 1960s  \cite{SLP1,SLP2,SLP3,SLP4,SLP5,S1} deals with the Fourier transform relating time (or spatial) dependent signals to their representation in frequency space,
namely

$$
Ff(k)=\int_{-\infty}^{\infty} e^{ikx}f(x) dx.
$$

  The crucial result is that for each pair of parameters $T, \Omega$, the integral operator of ``time-and-band limiting'' given by

$$
	U_{T,\Omega}f(x)=2\int_{-T}^T\frac{\sin \Omega(x-y)}{x-y}f(y) dy,
$$

commutes with the second order {\bf prolate spheroidal} differential operator

$$
R(x, \partial_x)= \partial_x (T^2-x^2)\partial_x -\Omega^2 x^2.
$$

	This result is fundamental in enabling one to compute with great accuracy the eigenfunctions of the integral operator, a very important issue in signal processing, geophysics and many other areas, see for example \cite{CG5,CG6,PS,SD,SDW}.
	For nontrivial connections with other areas see, for instance, \cite{CM,G9,FeldVa}. For more recent treatment of the numerical issues, see \cite{ORokX}.
For a fairly up-to-date account of applications see part V of the doctoral thesis of Yoan Le Henaff at Universit\'e de Rennes, France, pages 220--313. This includes extensive numerical work, \cite{LEH}.

\bigskip

 Since the key ingredient above, the exponential kernel $e^{i k x}$, satisfies a differential equation in $x$ and one in the dual variable $k$, one of us, \cite{G3, G4}, considered the situation one gets when one replaces a Fourier expansion by an expansion in terms of orthonormal polynomials $p_n(x)$ satisfying a ``bispectral property'', i.e. the existence not only of a three term recursion in $n$ but also the existence of a second order differential operator that would have all the $p_n(x)$ as its common eigenfunctions.

\bigskip

These so-called ``classical polynomials'' had been considered at least as far back as S. Bochner, \cite{Boch}.
\bigskip

Thus in \cite{G3} one looks at the Hermite, Laguerre and Jacobi cases, and for each pair of parameters $T,\Omega$ considers the integral operator of ``time-and-band limiting'' and in each of the three cases finds a differential operator that commutes with the integral one.

\bigskip

In the classical case of the Fourier transform the variables $x,k$ enter in a totally symmetric fashion. The case of orthogonal polynomials is not symmetric and one obtains as shown in the very nice papers by M. Perlstadt, see \cite{P1,P2}, that there is a very strong relation between a differential operator commuting with an integral one and tridiagonal matrices commuting with full matrices. This last situation is the one considered here.
This work has also been extended to some examples of ``exceptional orthogonal polynomials'', see \cite{GU2,CG7,CGZ}, and to matrix valued ones, see \cite{CG5,CG6,CGPZ,GPZ1,GPZ2,GPZ3}.

\bigskip

The search for these ``bispectral situations'' in the case when both spatial and frequency variables are continuous is given a detailed treatment in \cite{DG}.  A general method to connect bispectrality and the commutativity property
is given in \cite{GY} and was also developed further in \cite{CGYZ1,CGYZ2,CGYZ4}.

\bigskip

In this paper we replace the classical orthogonal polynomials of Bochner by the so-called symmetric classical orthogonal polynomials of Ch. Dunkl. See for instance \cite{Bel,ChGa,Dunk89}. These polynomials are common eigenfunctions of a second order {\bf Dunkl operator} defined in terms of Dunkl's first order differential-difference operator

	$$ T_\mu   = \frac{d}{dx} + (\mu/x)(I-R)$$

	where $R$ stands for the reflection operator $Rf(x)=f(-x)$ and $\mu$ is a real parameter.

It is shown in \cite{ChGa} that there are only two such families of polynomials, namely the ``generalized Hermite polynomials'' and the ``generalized Gegenbauer polynomials''. The first ones were already considered in \cite{Sz}.

\bigskip

One has that these polynomials satisfy

$$ (T_\mu^2 -2 x T_\mu) p_n(x)= \Lambda_n p_n(x)$$ in the Hermite case, and

$$((1-x^2) T_\mu^2 -2 (\alpha+1) x T_\mu) p_n(x)= \Lambda_n p_n(x)$$ in the Gegenbauer case, for an appropriate choice of $\Lambda_n$.

\bigskip

These polynomials also satisfy second order differential equations but the coefficients in the differential operator are $n$-dependent, ruling them out as ``classical orthogonal polynomials''; for instance, \cite{Sz} shows that the generalized Hermite polynomials $H_n^{\mu}$ satisfy
$$ x\,\frac{d^2}{dx^2}H_n^{\mu}(x) + 2( \mu -x^2)\frac{d}{dx}H_n^{\mu}(x)+ \left(2 n x-\frac{2\mu\theta_n}{x}\right) H_n^{\mu}(x)=0,$$
where $\theta_n=0$ if $n$ is even and $\theta_n=1$ if $n$ is odd.

\bigskip

These polynomials have appeared a long time ago \cite{Wi}, but they also play an important role in {\bf supersymmetric quantum mechanics} as formulated by E. Witten, see \cite{Witten1981}. See also \cite{CKS,PVZ,Qu,Qu23,TVZ,VZB}. Needless to say the full power of Dunkl operators is in the higher dimensional case. For a very interesting set of results in this situation, making contact with the bispectral problem, see \cite{Corr} and some of its references.

\bigskip

\section{The generalized Hermite polynomials}\label{GenHermite}

Let $\mu>-1/2$. We consider the family of generalized Hermite polynomials $H_n^{\mu}$, $n\geq 0$, given in terms of the confluent hypergeometric function ${}_1F_1$ by
\begin{equation}\label{defHmu}
	\begin{split}
		H_{2m}^{\mu}(x)&=(-1)^m\dfrac{(2m)!}{m!}\;{}_1F_1\!\left(-m;\mu+\tfrac12;x^2\right),\\
		H_{2m+1}^{\mu}(x)&=(-1)^m\dfrac{(2m+1)!}{m!\,\left(\mu+\frac12\right)}\;x\;{}_1F_1\!\left(-m;\mu+\tfrac32;x^2\right),
	\end{split}
\end{equation}
$m=0,1,2,\ldots$ (see for instance \cite{ChGa,Sz}). For $\mu=0$ one recovers, up to a normalization, the classical Hermite polynomials. Since ${}_1F_1(-m;\alpha+1;t)=\frac{m!}{(\alpha+1)_m}L_m^{(\alpha)}(t)$, these polynomials can also be written in terms of the classical Laguerre polynomials $L_m^{(\alpha)}$:
$$
H_{2m}^{\mu}(x)=(-1)^m\dfrac{(2m)!}{\left(\mu+\frac12\right)_m}\,L_m^{(\mu-\frac12)}(x^2),\qquad
H_{2m+1}^{\mu}(x)=(-1)^m\dfrac{(2m+1)!}{\left(\mu+\frac12\right)_{m+1}}\,x\,L_m^{(\mu+\frac12)}(x^2).
$$

\medskip

The generalized Hermite polynomials satisfy the orthogonality relation
\begin{equation}\label{orthHmu}
	\int_{-\infty}^{\infty}H_j^{\mu}(x)H_n^{\mu}(x)\,|x|^{2\mu}e^{-x^2}\,dx=h_n^{\mu}\,\delta_{n,j},
\end{equation}
where, writing $m=\lfloor n/2\rfloor$,
$$
h_n^{\mu}=\dfrac{(n!)^2\,\Gamma\!\left(\mu+\frac12\right)^2}{m!\,\Gamma\!\left(m+\mu+\frac12\right)}\quad (n \text{ even}),\qquad
h_n^{\mu}=\dfrac{(n!)^2\,\Gamma\!\left(\mu+\frac32\right)^2}{m!\,\Gamma\!\left(m+\mu+\frac32\right)\left(\mu+\frac12\right)^2}\quad (n \text{ odd}).
$$

\subsection{Bispectral operators}

Let us now consider the orthonormal sequence of generalized Hermite polynomials given by
$$\widetilde{H}_n^{\mu}=\dfrac{1}{\sqrt{h_n^{\mu}}}\,H_n^{\mu},\qquad n\geq 0 .$$

Let us define
\begin{equation}\label{gammas}
	\gamma_n=\begin{cases} n, & n \text{ even},\\ n+2\mu, & n \text{ odd}.\end{cases}
\end{equation}

The polynomials $\widetilde{H}_n^{\mu}$ satisfy the three term recurrence relation
\begin{equation}\label{recurrHmu}
	a_n\widetilde{H}_{n-1}^{\mu}(x)+a_{n+1}\widetilde{H}_{n+1}^{\mu}(x)=x\,\widetilde{H}_n^{\mu}(x),\qquad a_n=\sqrt{\dfrac{\gamma_n}{2}},\quad n\geq 0,
\end{equation}
with $\widetilde{H}_{-1}^{\mu}=0$.

\medskip

On the other hand, in terms of the Dunkl operator $T_\mu f(x)=f'(x)+\dfrac{\mu}{x}\left(f(x)-f(-x)\right)$, one has
$$
\left(T_\mu^2-2x\,T_\mu\right)\widetilde{H}_n^{\mu}=-2\gamma_n\,\widetilde{H}_n^{\mu},\qquad n\geq 0 .
$$

\subsection{Time-and-band limiting}\label{TBLHermite}

Consider the $N\times N$ matrix $M$ of truncated inner products depending on a real parameter $T>0$, whose entries are given by
\begin{equation}\label{Mmu}
	M_{m,n}=\int_{-T}^{T}\widetilde{H}_m^{\mu}(x)\widetilde{H}_n^{\mu}(x)\,|x|^{2\mu}e^{-x^2}dx,\quad m,n=0,1,\ldots, N-1.
\end{equation}

\bigskip

The issue of finding a commuting banded matrix has been considered by us in previous publications, see \cite{CG5, CG7, GPZ1, GPZ3, GVZ}.
In these papers we conjectured the general expression for the commuting matrix by extrapolation from several special
choices of $(T,N)$ and in many cases we eventually found a proof of the property in question.
In the present case we found ourselves the commuting matrix for several choices of the parameters $(T,N)$ and then
we asked Claude for assistance. The general formulas here were thus obtained with the help of Claude.

\bigskip

For fixed values of the parameters $N$, $\mu$ and $T$ one looks for a ``narrow banded'' commuting matrix $\mathscr{L}$. One finds that there exists a {\it pentadiagonal} symmetric matrix $\mathscr{L}$ commuting with $M$, which in fact has nonzero entries only on the main diagonal and on the second diagonals. Its entries are given explicitly, for every $N$, $\mu$ and $T$, by
\begin{equation}\label{Lmu}
	\begin{split}
		\mathscr{L}_{n,n}&=k_n\,(2n+2\mu+1)+n\,(T^2+1),\qquad n=0,1,\ldots,N-1,\\
		\mathscr{L}_{n,n+2}&=\mathscr{L}_{n+2,n}=k_n\,\sqrt{\gamma_{n+1}\,\gamma_{n+2}}\;=\;2\,k_n\,a_{n+1}a_{n+2},\qquad n=0,1,\ldots,N-3,
	\end{split}
\end{equation}
and $\mathscr{L}_{m,n}=0$ otherwise, where $k_n=\left\lfloor \dfrac{N-1-n}{2}\right\rfloor$ and $\gamma_n$, $a_n$ are given in (\ref{gammas}) and (\ref{recurrHmu}). Notice that $2n+2\mu+1=\gamma_n+\gamma_{n+1}=2(a_n^2+a_{n+1}^2)$.

\bigskip

For instance, for $N=6$ we have the following expression for the commuting matrix $\mathscr{L}$, valid for arbitrary values of $\mu$ and $T$:
$$
\mathscr{L}=\left(
\begin{array}{cccccc}
	4\mu+2 & 0 & 2\sqrt{2(2\mu+1)} & 0 & 0 & 0 \\
	0 & T^2+4\mu+7 & 0 & 2\sqrt{2(2\mu+3)} & 0 & 0 \\
	2\sqrt{2(2\mu+1)} & 0 & 2T^2+2\mu+7 & 0 & 2\sqrt{2\mu+3} & 0 \\
	0 & 2\sqrt{2(2\mu+3)} & 0 & 3T^2+2\mu+10 & 0 & 2\sqrt{2\mu+5} \\
	0 & 0 & 2\sqrt{2\mu+3} & 0 & 4(T^2+1) & 0 \\
	0 & 0 & 0 & 2\sqrt{2\mu+5} & 0 & 5(T^2+1)
\end{array}
\right).
$$

\bigskip
If one could prove that $\mathscr{L}$ has simple spectrum we would be done. This is generically true but it is not true for some values of the parameters. For instance if $N=7,\mu=3$ there are two values of $T$, where the resultant of its characteristic polynomial and its derivative vanishes simultaneously. 
Fortunately this condition of simple spectrum is not really needed as a careful analysis of $\mathscr{L}$ given below reveals.

 \bigskip

Since $\mathscr{L}_{m,n}=0$ and $M_{m,n}=0$ whenever $m+n$ is odd, both split into an even block (indices $0,2,4,\dots$) and an odd block (indices $1,3,5,\dots$):
\[
\mathscr{L}=\mathscr{L}_{\mathrm{even}}\oplus \mathscr{L}_{\mathrm{odd}},\qquad M=M_{\mathrm{even}}\oplus M_{\mathrm{odd}},
\] and the characteristic polynomial of $\mathscr{L}$ factorizes as
$\chi_{\mathscr{L}}=\chi_{\mathrm{even}}\,\chi_{\mathrm{odd}}$.

\bigskip
The blocks give rise to tridiagonal matrices with non-zero off-diagonals and thus simple spectrum.

\bigskip
This and the previous observations lead to

\begin{thm} \label{thmeig} 
	For all $N\ge1$, $T>0$, $\mu>-\tfrac12$, every eigenvector of $\mathscr{L}_{\mathrm{even}}$ (extended by zero to the odd indices) and every eigenvector of $\mathscr{L}_{\mathrm{odd}}$ (extended by zero to the even indices) is an eigenvector of $M$. These $N$ vectors form an orthonormal eigenbasis of $M$.
\end{thm}
\begin{proof}
	$M$ is block-diagonal by parity and $\mathscr{L}M=M\mathscr{L}$, so $\mathscr{L}_{\mathrm{even}}M_{\mathrm{even}}=M_{\mathrm{even}}\mathscr{L}_{\mathrm{even}}$, and likewise for the odd blocks. A matrix commuting with a matrix of simple spectrum preserves each of its one-dimensional eigenspaces, so each eigenvector of $\mathscr{L}_{\mathrm{even}}$ is an eigenvector of $M_{\mathrm{even}}$, and likewise for the odd block.
\end{proof}

To sum up: the eigenvectors of $M$ are obtained from $\mathscr{L}$ for every $N$, $T$, $\mu$, whether or not the full spectrum of $\mathscr{L}$ is simple.

\subsection{The commuting Dunkl operator}\label{proofcomm}

Claude gives a proof of the commutativity that parallels the arguments given in the classical cases in \cite{G3,P1,P2}. We summarize the main points.

As in the classical situation, the key point is that $\mathscr{L}$ is the matrix, in the basis $\{\widetilde{H}^{\mu}_n\}$, of a second order operator which is symmetric on $(-T,T)$ and preserves the space of polynomials of degree less than $N$. In the present case this operator is not a differential operator but a Dunkl operator, involving the reflection $Rf(x)=f(-x)$.

\medskip

Let $\Pi_0=\frac12(I+R)$ and $\Pi_1=\frac12(I-R)$ be the projections onto even and odd functions, and let $N_0$ (respectively $N_1$) be the largest even (respectively odd) integer not greater than $N-1$. Notice that $k_n=\frac{N_\epsilon-n}{2}$ if $n\equiv\epsilon \pmod 2$.

The crucial point in the argument is that the operator
\begin{equation}\label{opDmu}
	\mathcal{D}=-\frac12\, e^{x^2}\,T_\mu\, e^{-x^2}\left(T^2-x^2\right)T_\mu\;+\;x^2\left(\gamma_{N_0}\,\Pi_0+\gamma_{N_1}\,\Pi_1\right)\;-\;2\mu\left(T^2+1\right)\Pi_1 ,
\end{equation} has the following property:

\begin{thm}\label{thmcomm}
	For every $N\ge1$, $\mu>-1/2$ and $T>0$ one has
	\begin{equation}\label{Dmatrix}
		\mathcal{D}\,\widetilde{H}^{\mu}_n=\sum_{m=0}^{N-1}\mathscr{L}_{m,n}\,\widetilde{H}^{\mu}_m,\qquad n=0,1,\ldots,N-1,
	\end{equation}
	where $\mathscr{L}$ is the matrix given in (\ref{Lmu}). With a bit of extra work one can get the consequence that $\mathscr{L}M=M\mathscr{L}$.
\end{thm}

\begin{remark}
	The same argument shows that $\mathcal{D}$ commutes with the integral operator $\mathcal{S}$ acting on functions on $(-T,T)$ by means of the kernel $k(x,y)=\sum_{j=0}^{N-1}\widetilde{H}^{\mu}_j(x)\widetilde{H}_j^{\mu}(y)$. Indeed, by (\ref{Dmatrix}) and the symmetry of $\mathscr{L}$, one has $\mathcal{D}_x k(x,y)=\sum_{j,m=0}^{N-1}\mathscr{L}_{m,j}\,\widetilde{H}_m^{\mu}(x)\widetilde{H}_j^{\mu}(y)=\mathcal{D}_y k(x,y)$, and then one argues that $\mathcal{D}\mathcal{S}f=\mathcal{S}\mathcal{D}f$. Here $\mathcal{D}$ plays the role of the differential operator found by Slepian, Landau and Pollak, or in references \cite{G3,P1,P2}.
\end{remark}

\section{The generalized Gegenbauer polynomials}\label{GenGegenbauer}

Let $\alpha>-1$ and $\mu>-1/2$. Following \cite{ChGa} (see also \cite{Bel}), we consider the family of generalized Gegenbauer polynomials $S_n=S_n^{(\alpha,\mu-\frac12)}$, $n\geq0$, given explicitly by
\begin{equation}\label{defSmu}
	S_n(x)=\sum_{k=0}^{\lfloor n/2\rfloor}\frac{(-1)^k\left(\alpha+\mu+\frac12\right)_{n-k}}{k!\;c_\mu(n-2k)}\,(2x)^{n-2k},
\end{equation}
where $c_\mu(2m)=2^{2m}m!\left(\mu+\tfrac12\right)_m$ and $c_\mu(2m+1)=2^{2m+1}m!\left(\mu+\tfrac12\right)_{m+1}$.
For $\mu=0$ one recovers, up to a normalization, the classical Gegenbauer polynomials $C_n^{(\alpha+\frac12)}$. As in the Hermite case, the even and odd polynomials can be written in terms of classical orthogonal polynomials in the variable $x^2$, in this case the Jacobi polynomials $P_m^{(a,b)}$:
$$
S_{2m}(x)=\frac{\left(\alpha+\mu+\frac12\right)_m}{\left(\mu+\frac12\right)_m}\,P_m^{(\alpha,\mu-\frac12)}(2x^2-1),\qquad
S_{2m+1}(x)=\frac{\left(\alpha+\mu+\frac12\right)_{m+1}}{\left(\mu+\frac12\right)_{m+1}}\,x\,P_m^{(\alpha,\mu+\frac12)}(2x^2-1).
$$

\medskip

The generalized Gegenbauer polynomials satisfy the orthogonality relation
\begin{equation}\label{orthSmu}
	\int_{-1}^{1}S_j(x)S_n(x)\,|x|^{2\mu}(1-x^2)^{\alpha}\,dx=h_n\,\delta_{n,j},
\end{equation}
where, writing $m=\lfloor n/2\rfloor$ and $\theta_n=0$ if $n$ is even, $\theta_n=1$ if $n$ is odd,
$$
h_n=\frac{\Gamma\!\left(\mu+\frac12\right)^2\,\Gamma\!\left(m+\alpha+\mu+\frac12+\theta_n\right)\Gamma(m+\alpha+1)}{\Gamma\!\left(\alpha+\mu+\frac12\right)^2\,\Gamma\!\left(m+\mu+\frac12+\theta_n\right)\,m!\,\left(2m+\alpha+\mu+\frac12+\theta_n\right)} .
$$

\subsection{Bispectral operators}

Let us consider the orthonormal sequence of generalized Gegenbauer polynomials
$$\widetilde{S}_n=\frac{1}{\sqrt{h_n}}\,S_n,\qquad n\geq0.$$
If $\alpha+\mu+\frac12<0$ the leading coefficient of $S_n$, $n\ge1$, is negative, and in that case we replace $S_n$ by $-S_n$ for $n\geq1$, so that all the $\widetilde{S}_n$ have positive leading coefficient. If $\alpha+\mu+\frac12=0$ the polynomials $S_n$, $n\ge1$, given by (\ref{defSmu}) vanish identically; in this case the orthonormal polynomials $\widetilde{S}_n$ are defined by continuity in the parameters, and all the formulas below remain valid in the same sense.

With $\gamma_n$ as in (\ref{gammas}), the polynomials $\widetilde{S}_n$ satisfy the three term recurrence relation
\begin{equation}\label{recurrSmu}
	\begin{split}
	&a_n\widetilde{S}_{n-1}(x)+a_{n+1}\widetilde{S}_{n+1}(x)=x\,\widetilde{S}_n(x),\\
	&a_n=\sqrt{\frac{\gamma_n\left(\gamma_n+2\alpha\right)}{(2n+2\alpha+2\mu-1)(2n+2\alpha+2\mu+1)}},\quad n\geq0,
	\end{split}
\end{equation} with $\widetilde{S}_{-1}=0$, i.e. $J\Psi(x)=x\Psi(x)$, where $\Psi(x)=(\widetilde{S}_0(x),\widetilde{S}_1(x),\widetilde{S}_2(x),\ldots)^t$ and $J$ is the semi-infinite symmetric tridiagonal matrix with zero diagonal and off-diagonal entries $a_1,a_2,a_3,\ldots$.

\medskip

On the other hand, in terms of the Dunkl operator $T_\mu$, one has
\begin{equation}\label{DunklSmu}
	\left((1-x^2)\,T_\mu^2-2(\alpha+1)\,x\,T_\mu\right)\widetilde{S}_n=\Lambda_n\,\widetilde{S}_n,\qquad
	\Lambda_n=-\gamma_n\left(\gamma_{n-1}+2\alpha+2\right),\qquad n\geq0,
\end{equation}
that is, $\Lambda_n=-n(n+2\alpha+2\mu+1)$ if $n$ is even and $\Lambda_n=-(n+2\mu)(n+2\alpha+1)$ if $n$ is odd.

\subsection{Time-and-band limiting}\label{TBLGegenbauer}

For a real parameter $0<T<1$ we consider the $N\times N$ matrix $M$ of truncated inner products, whose entries are given by
\begin{equation}\label{MSmu}
	M_{m,n}=\int_{-T}^{T}\widetilde{S}_m(x)\widetilde{S}_n(x)\,|x|^{2\mu}(1-x^2)^{\alpha}\,dx,\qquad m,n=0,1,\ldots,N-1.
\end{equation}
As in the Hermite case, $M_{m,n}=0$ whenever $m+n$ is odd.

For $T=1$ one gets the identity matrix.

\bigskip

As in the Hermite case, for any $N$, $\alpha$, $\mu$ and $T$ there exists a {\it pentadiagonal} symmetric matrix $\mathscr{L}$ commuting with $M$, with nonzero entries only on the main diagonal and on the second diagonals. With $k_n=\left\lfloor\frac{N-1-n}{2}\right\rfloor$ and
$$
\lambda_n=4k_n\left(n+k_n+\alpha+\mu+\tfrac32\right),
$$
its entries are given by
\begin{equation}\label{LSmu}
	\begin{split}
		\mathscr{L}_{n,n}&=\lambda_n\left(a_n^2+a_{n+1}^2\right)+2\left(2n+2\alpha+2\mu+1\right)a_n^2\\
		&\quad+n\left(n+2\alpha+2\mu+1\right)T^2-2\mu\,\theta_n,\qquad n=0,1,\ldots,N-1,\\
		\mathscr{L}_{n,n+2}&=\mathscr{L}_{n+2,n}=\lambda_n\,a_{n+1}a_{n+2},\qquad n=0,1,\ldots,N-3,
	\end{split}
\end{equation}
and $\mathscr{L}_{m,n}=0$ otherwise, where $a_n$ is given in (\ref{recurrSmu}). Notice that $2\left(2n+2\alpha+2\mu+1\right)a_n^2=\dfrac{2\gamma_n(\gamma_n+2\alpha)}{2n+2\alpha+2\mu-1}$. The analogy with (\ref{Lmu}) becomes clear if one writes the Hermite matrix as $\mathscr{L}_{n,n}=2k_n(a_n^2+a_{n+1}^2)+2a_n^2+nT^2+n-\gamma_n$, $\mathscr{L}_{n,n+2}=2k_na_{n+1}a_{n+2}$, and notices that $n-\gamma_n=-2\mu\theta_n$.

\bigskip

The expressions of the entries of $\mathscr{L}$ for general $\alpha$ and $\mu$ are rather long. For instance, for $N=6$ and $\alpha=\mu=1$ we have the following expression, valid for arbitrary values of $T$:
$$
\mathscr{L}=\left(
\begin{array}{cccccc}
	\frac{132}{7} & 0 & \frac{88\sqrt{6}}{21} & 0 & 0 & 0 \\[1mm]
	0 & 6T^2+\frac{296}{9} & 0 & \frac{104\sqrt{110}}{99} & 0 & 0 \\[1mm]
	\frac{88\sqrt{6}}{21} & 0 & 14T^2+\frac{1138}{77} & 0 & \frac{4\sqrt{2730}}{33} & 0 \\[1mm]
	0 & \frac{104\sqrt{110}}{99} & 0 & 24T^2+\frac{2506}{117} & 0 & \frac{36\sqrt{770}}{143} \\[1mm]
	0 & 0 & \frac{4\sqrt{2730}}{33} & 0 & 36T^2+\frac{48}{11} & 0 \\[1mm]
	0 & 0 & 0 & \frac{36\sqrt{770}}{143} & 0 & 50T^2+\frac{100}{13}
\end{array}
\right).
$$

\bigskip


Using an argument similar to that given at the end of the previous section, one can show that the eigenvectors of $M$ are obtained from $\mathscr{L}$ for every $N$, $T$, $\mu$,  $\alpha$, whether or not the full spectrum of $\mathscr{L}$ is simple.

 As in the preceding section, the existence of a commuting banded operator was first established by us for several special choices of the parameters  and the general formulas were obtained with the help of Claude.

\subsection{The commuting Dunkl operator}\label{proofcommG}

There are very few differences between the proof in this case and in
the Hermite one. $\mathscr{L}$ is the matrix, in the basis $\{\widetilde{S}_n\}$, of a second order Dunkl operator which is symmetric on $(-T,T)$ and preserves the space of polynomials of degree less than $N$.

One keeps the notation $\Pi_0$, $\Pi_1$, $N_0$, $N_1$ of Subsection \ref{proofcomm}, and puts
$$
\rho_n=\gamma_n\left(\gamma_{n-1}+2\alpha+4\right)=2\gamma_n-\Lambda_n .
$$

One defines the operator
\begin{equation}\label{opDSmu}
	\begin{split}
	\mathcal{D}&=-(1-x^2)^{-\alpha}\,T_\mu\,(1-x^2)^{\alpha+1}\left(T^2-x^2\right)T_\mu\\
	&\quad+\;x^2\left(\rho_{N_0}\,\Pi_0+\rho_{N_1}\,\Pi_1\right)\;-\;2\mu\left((2\alpha+1)T^2+1\right)\Pi_1 ,
	\end{split}
\end{equation}
which, as in the Hermite case, has the following property:

\begin{thm}\label{thmcommG}
	For every $N\ge1$, $\alpha>-1$, $\mu>-1/2$ and $0<T<1$ one has
	\begin{equation}\label{DSmatrix}
		\mathcal{D}\,\widetilde{S}_n=\sum_{m=0}^{N-1}\mathscr{L}_{m,n}\,\widetilde{S}_m,\qquad n=0,1,\ldots,N-1,
	\end{equation}
	where $\mathscr{L}$ is the matrix given in (\ref{LSmu}). Once again with some extra labor one gets the consequence $\mathscr{L}M=M\mathscr{L}$.
\end{thm}

\begin{remark}
	As in the Hermite case and in \cite{G3,P1,P2}, $\mathcal{D}$ also commutes with the integral operator on $(-T,T)$ with kernel $k(x,y)=\sum_{j=0}^{N-1}\widetilde{S}_j(x)\widetilde{S}_j(y)$.
\end{remark}

\section{The benefit of having a commuting local matrix}\label{numeric}

The previous sections have exhibited, for generalized Hermite and Gegenbauer polynomials respectively, a pair of matrices. The first one, a full matrix $M=M_{T,N}$, is obtained by forming the inner products of the normalized orthogonal polynomials over a restricted range in physical space, thus implementing ``time limiting'' with parameter $T$. In the resulting $N \times N$ matrix the parameter $N$ implements ``band limiting''. In each case the second matrix, denoted by $\mathscr{L}=\mathscr{L}_{T,N}$, is a narrow banded one that commutes with the first matrix, and whose even and odd blocks have simple spectrum. From a numerical point of view the entire purpose of the search for this second matrix is that it reduces the problem of computing the eigenvectors of the first one, a seriously ill-conditioned one, to a very well conditioned one.

\bigskip

The eigenvectors of $M_{T,N}$ are of paramount importance since they give the singular vectors of the signal processing problem at hand, as described in Section \ref{intro}.

\bigskip

Following the approach in \cite[Section 6]{CG5}, \cite[Section 7]{CG7}  and \cite[Section 5]{CFF}, we display below the results of some relatively small size numerical computations that illustrate the problem of computing the eigenvectors of a full matrix, such as $M_{T,N}$, some of whose eigenvalues are very close together. In each case, we give the eigenvalues of both the full matrix  $M_{T,N}$ and those of the narrow banded matrix $\mathscr{L}_{T,N}$.  It should be clear that in the case when $N$ is large the problems indicated below get to be much worse. Our point is that they already appear for relatively small values of $N$.

\bigskip

We use the QR algorithm as implemented in LAPACK. In each of the two situations  we will denote	by $X_M$ the matrix of eigenvectors  of $M_{T,N}$ (normalized and given as columns of $X_M$). We will denote by $Y_L$ the matrix of eigenvectors  of $\mathscr{L}_{T,N}$ (normalized and given as columns of $Y_L$).

\bigskip

In theory the eigenvectors of $M_{T,N}$  should agree (up to order and signs) with those of $\mathscr{L}_{T,N}$. If we compute the matrix
of inner products given by $$Y_L^T X_M$$ we expect to have the identity matrix up to some permutation and possibly some signs due to the normalization of the eigenvectors which are the columns of $X_M$ and $Y_L$.

\subsection{Generalized Hermite}

The choice of parameters is $N=12$, $\mu=4$ and $T=7$.

\bigskip

The eigenvalues of $M_{T,N}$, as computed in double precision, are
$$
\begin{array}{c}
	1.0000000000000009,\; 1.0000000000000004,\; 1.0000000000000002,\; 0.9999999999999999,\\
	0.9999999999999999,\; 0.9999999999999994,\; 0.9999999999999994,\; 0.9999999999999929,\\
	0.9999999999802371,\; 0.9999999998902893,\; 0.9999995186658811,\; 0.9999980178732023,
\end{array}
$$
\noindent
and those of $\mathscr{L}_{T,N}$ are
$$
\begin{array}{c}
	552.636,\;502.290,\;481.022,\;428.447,\;398.303,\;343.875,\\
	306.489,\;250.413,\;206.734,\;149.150,\;99.8169,\;40.8251.
\end{array}
$$

\bigskip

In order to make visible the numerical errors, we display the matrix whose entries are the absolute values of the entries of $Y_L^T X_M$, that is $\left|Y_L^T X_M\right|$. Here the columns of $X_M$ are ordered according to the eigenvalues of $M_{T,N}$ in decreasing order, and those of $Y_L$ according to the eigenvalues of $\mathscr{L}_{T,N}$ in decreasing order. For the sake of space we split it in two blocks. The first seven columns of $\left|Y_L^T X_M\right|$ are

\begin{scriptsize}
	$$
	\left(
	\begin{array}{ccccccc}
		5.10\times 10^{-14} & 3.43\times 10^{-14} & 7.22\times 10^{-15} & 5.95\times 10^{-17} & 1.29\times 10^{-16} & 2.28\times 10^{-15} & 2.80\times 10^{-14} \\
		5.90\times 10^{-12} & 4.90\times 10^{-12} & 1.49\times 10^{-11} & 8.17\times 10^{-14} & 5.96\times 10^{-13} & 1.06\times 10^{-12} & 2.64\times 10^{-11} \\
		8.73\times 10^{-9} & 2.61\times 10^{-8} & 2.41\times 10^{-8} & 2.76\times 10^{-13} & 2.01\times 10^{-10} & 6.42\times 10^{-10} & 6.74\times 10^{-9} \\
		4.90\times 10^{-6} & 1.59\times 10^{-5} & 9.14\times 10^{-6} & 1.82\times 10^{-8} & 2.46\times 10^{-7} & 2.29\times 10^{-7} & 5.23\times 10^{-6} \\
		2.93\times 10^{-3} & 1.25\times 10^{-3} & 6.92\times 10^{-4} & 6.72\times 10^{-7} & 2.16\times 10^{-5} & 4.89\times 10^{-5} & 1.04\times 10^{-3} \\
		7.83\times 10^{-4} & 0.249519 & 0.231040 & 1.25\times 10^{-3} & 1.06\times 10^{-3} & 0.038341 & 0.939620 \\
		0.282708 & 0.822301 & 0.483360 & 3.32\times 10^{-4} & 9.76\times 10^{-4} & 0.013928 & 0.100318 \\
		0.322779 & 0.481870 & 0.697184 & 0.030337 & 0.274407 & 0.098256 & 0.302780 \\
		0.878144 & 0.086437 & 0.391164 & 6.60\times 10^{-3} & 0.010339 & 0.247306 & 0.084052 \\
		0.075197 & 0.139625 & 0.209989 & 0.203216 & 0.938771 & 0.012907 & 0.089452 \\
		0.197609 & 0.047322 & 0.171630 & 0.033554 & 0.037055 & 0.962555 & 0.015571 \\
		6.37\times 10^{-3} & 0.012454 & 0.018889 & 0.978065 & 0.204765 & 0.029016 & 8.06\times 10^{-3} \\
	\end{array}
	\right)
	$$
\end{scriptsize}

\noindent
and the last five columns are

\begin{scriptsize}
	$$
	\left(
	\begin{array}{ccccc}
		1.03\times 10^{-13} & 2.99\times 10^{-13} & 3.55\times 10^{-12} & 2.24\times 10^{-12} & 1.000000 \\
		1.24\times 10^{-10} & 1.18\times 10^{-10} & 5.27\times 10^{-11} & 1.000000 & 2.24\times 10^{-12} \\
		3.20\times 10^{-7} & 4.97\times 10^{-8} & 1.000000 & 5.27\times 10^{-11} & 3.55\times 10^{-12} \\
		2.21\times 10^{-6} & 1.000000 & 4.97\times 10^{-8} & 1.18\times 10^{-10} & 2.99\times 10^{-13} \\
		0.999994 & 2.25\times 10^{-6} & 3.20\times 10^{-7} & 1.23\times 10^{-10} & 1.04\times 10^{-13} \\
		1.45\times 10^{-3} & 6.76\times 10^{-6} & 4.97\times 10^{-9} & 2.72\times 10^{-11} & 3.65\times 10^{-14} \\
		1.42\times 10^{-3} & 1.83\times 10^{-5} & 3.17\times 10^{-8} & 1.20\times 10^{-11} & 3.62\times 10^{-14} \\
		1.65\times 10^{-4} & 1.97\times 10^{-6} & 4.27\times 10^{-10} & 1.93\times 10^{-12} & 3.14\times 10^{-15} \\
		2.66\times 10^{-3} & 1.15\times 10^{-6} & 1.81\times 10^{-8} & 9.80\times 10^{-12} & 4.26\times 10^{-14} \\
		1.45\times 10^{-5} & 3.00\times 10^{-7} & 2.46\times 10^{-11} & 2.07\times 10^{-13} & 1.81\times 10^{-16} \\
		6.94\times 10^{-4} & 2.80\times 10^{-7} & 3.69\times 10^{-9} & 2.13\times 10^{-12} & 1.02\times 10^{-14} \\
		7.39\times 10^{-7} & 2.23\times 10^{-8} & 8.29\times 10^{-13} & 4.33\times 10^{-14} & 6.09\times 10^{-18} \\
	\end{array}
	\right)
	$$
\end{scriptsize}

\bigskip

The last five columns, which correspond to the eigenvalues of $M_{T,N}$ that can still be separated in double precision, give (up to the order) the identity matrix, with the remaining entries close to zero. On the contrary, the first seven columns correspond to eigenvalues of $M_{T,N}$ that are equal to $1$ to machine precision, and the computed eigenvectors of $M_{T,N}$ are just arbitrary orthonormal combinations of the true ones. Indeed, computing with high precision one obtains that the eigenvalues of $M_{T,N}$ are $1-\varepsilon_k$, $k=1,\ldots,12$, with
$$
\begin{array}{c}
	\varepsilon_k=1.98\times 10^{-6},\; 4.81\times 10^{-7},\; 1.10\times 10^{-10},\; 1.98\times 10^{-11},\; 7.23\times 10^{-15},\; 1.00\times 10^{-15},\\ 3.77\times 10^{-19},\; 4.08\times 10^{-20},\; 1.22\times 10^{-23},\; 1.02\times 10^{-24},\; 1.73\times 10^{-28},\; 1.09\times 10^{-29},
\end{array}
$$
so that seven of them differ from $1$ by an amount of order $10^{-15}$ or less, that is, at the level of the rounding errors of double precision arithmetic, and they cannot be told apart. The eigenvalues of $\mathscr{L}_{T,N}$, instead, are well separated and its eigenvectors are computed without any difficulty.

\subsection{Generalized Gegenbauer}

The choice of parameters is $N=12$, $\alpha=1$, $\mu=1$ and $T=1/16$.

\bigskip

In this case $T$ is small, and the eigenvalues of $M_{T,N}$ accumulate at $0$ instead of at $1$. This will have an effect later on: the role of the first columns in the Hermite case will be taken by the last columns in the Gegenbauer case.  As computed in double precision, they are
$$
\begin{array}{c}
	0.03869726063876869,\; 0.0015466977718477978,\; 2.571427024967364\times 10^{-6},\\
	2.0304594441879703\times 10^{-8},\; 1.0488849635883766\times 10^{-11},\; 3.913604968931161\times 10^{-14},\\
	8.079300388352738\times 10^{-18},\; -1.156892006873053\times 10^{-18},\; 4.9259485961842345\times 10^{-19},\\
	-4.050965345287995\times 10^{-19},\; -6.690456697819793\times 10^{-20},\; 5.881582269125799\times 10^{-20},
\end{array}
$$
\noindent
and those of $\mathscr{L}_{T,N}$ are
$$
\begin{array}{c}
	154.091,\;130.075,\;108.184,\;88.1570,\;70.2628,\;54.2242,\\
	40.3309,\;28.2767,\;18.3996,\;10.3178,\;4.55163,\;0.395210.
\end{array}
$$

\bigskip

As before, we display the matrix $\left|Y_L^T X_M\right|$. Here the columns of $X_M$ are ordered according to the computed eigenvalues of $M_{T,N}$, in decreasing order of their absolute values (as listed above), and those of $Y_L$ according to the eigenvalues of $\mathscr{L}_{T,N}$ in decreasing order. The first six columns of $\left|Y_L^T X_M\right|$ are

\begin{scriptsize}
	$$
	\left(
	\begin{array}{cccccc}
		1.58\times 10^{-17} & 2.16\times 10^{-17} & 7.26\times 10^{-15} & 8.46\times 10^{-12} & 8.43\times 10^{-10} & 1.89\times 10^{-6} \\
		4.87\times 10^{-18} & 2.18\times 10^{-18} & 2.36\times 10^{-14} & 4.78\times 10^{-12} & 4.74\times 10^{-10} & 5.33\times 10^{-6} \\
		4.06\times 10^{-17} & 2.79\times 10^{-16} & 6.22\times 10^{-15} & 3.60\times 10^{-12} & 2.02\times 10^{-8} & 5.95\times 10^{-6} \\
		1.12\times 10^{-16} & 1.68\times 10^{-17} & 5.10\times 10^{-14} & 1.48\times 10^{-11} & 3.50\times 10^{-8} & 9.28\times 10^{-6} \\
		8.48\times 10^{-18} & 5.71\times 10^{-16} & 1.01\times 10^{-13} & 1.78\times 10^{-11} & 3.04\times 10^{-8} & 1.37\times 10^{-6} \\
		3.41\times 10^{-16} & 6.97\times 10^{-17} & 1.63\times 10^{-13} & 1.34\times 10^{-11} & 4.30\times 10^{-8} & 5.88\times 10^{-6} \\
		3.26\times 10^{-16} & 7.36\times 10^{-16} & 2.32\times 10^{-13} & 2.82\times 10^{-11} & 3.57\times 10^{-8} & 1.000000 \\
		1.45\times 10^{-17} & 1.52\times 10^{-16} & 3.37\times 10^{-13} & 2.55\times 10^{-12} & 1.000000 & 3.57\times 10^{-8} \\
		4.79\times 10^{-17} & 6.72\times 10^{-16} & 4.64\times 10^{-14} & 1.000000 & 2.55\times 10^{-12} & 2.82\times 10^{-11} \\
		2.06\times 10^{-16} & 2.60\times 10^{-16} & 1.000000 & 4.63\times 10^{-14} & 3.37\times 10^{-13} & 2.32\times 10^{-13} \\
		7.01\times 10^{-17} & 1.000000 & 2.30\times 10^{-16} & 8.81\times 10^{-16} & 1.95\times 10^{-16} & 8.70\times 10^{-16} \\
		1.000000 & 1.44\times 10^{-17} & 4.86\times 10^{-17} & 9.26\times 10^{-18} & 2.47\times 10^{-16} & 7.20\times 10^{-17} \\
	\end{array}
	\right)
	$$
\end{scriptsize}

\noindent
and the last six columns are

\begin{scriptsize}
	$$
	\left(
	\begin{array}{cccccc}
		5.10\times 10^{-3} & 0.329219 & 0.352524 & 0.343996 & 0.412034 & 0.692250 \\
		0.251383 & 0.056429 & 0.843861 & 0.065916 & 0.454893 & 0.101234 \\
		0.107154 & 0.576653 & 0.155952 & 0.332972 & 0.260112 & 0.673156 \\
		0.147256 & 0.364927 & 0.274449 & 0.426098 & 0.741307 & 0.196790 \\
		0.024185 & 0.639294 & 0.044936 & 0.752172 & 0.070174 & 0.134212 \\
		0.950278 & 0.118463 & 0.248924 & 0.138309 & 0.034367 & 0.025763 \\
		4.90\times 10^{-6} & 6.96\times 10^{-6} & 2.14\times 10^{-6} & 8.08\times 10^{-6} & 6.68\times 10^{-6} & 1.75\times 10^{-6} \\
		3.87\times 10^{-8} & 2.58\times 10^{-8} & 2.41\times 10^{-8} & 8.48\times 10^{-9} & 3.17\times 10^{-8} & 2.40\times 10^{-8} \\
		1.10\times 10^{-11} & 1.77\times 10^{-13} & 9.92\times 10^{-13} & 2.54\times 10^{-11} & 7.01\times 10^{-12} & 3.12\times 10^{-12} \\
		1.65\times 10^{-13} & 5.93\times 10^{-14} & 1.50\times 10^{-15} & 8.00\times 10^{-14} & 4.42\times 10^{-14} & 2.93\times 10^{-14} \\
		8.36\times 10^{-17} & 3.38\times 10^{-16} & 8.25\times 10^{-18} & 3.28\times 10^{-16} & 1.96\times 10^{-17} & 5.92\times 10^{-17} \\
		6.54\times 10^{-16} & 5.89\times 10^{-17} & 3.09\times 10^{-16} & 2.45\times 10^{-16} & 4.53\times 10^{-18} & 3.85\times 10^{-17} \\
	\end{array}
	\right)
	$$
\end{scriptsize}

\bigskip

The first six columns, which correspond to the eigenvalues of $M_{T,N}$ that can still be resolved in double precision, give (up to the order) the identity matrix, with the remaining entries close to zero. On the contrary, the last six columns correspond to eigenvalues of $M_{T,N}$ that are equal to $0$ to machine precision, and the computed eigenvectors are just arbitrary orthonormal combinations of the true ones. Indeed, computing with high precision one obtains that the eigenvalues of $M_{T,N}$ are
$$
\begin{array}{c}
	3.87\times 10^{-2},\; 1.55\times 10^{-3},\; 2.57\times 10^{-6},\; 2.03\times 10^{-8},\; 1.05\times 10^{-11},\; 3.91\times 10^{-14},\\
	6.98\times 10^{-18},\; 1.64\times 10^{-20},\; 9.48\times 10^{-25},\; 1.62\times 10^{-27},\; 2.18\times 10^{-32},\; 2.94\times 10^{-35},
\end{array}
$$
so that the last six of them are smaller than $10^{-17}$, that is, below the level of the rounding errors of double precision arithmetic relative to the largest one; the values computed in double precision for them are meaningless (three of them even have the wrong sign). The eigenvalues of $\mathscr{L}_{T,N}$, instead, are well separated and its eigenvectors are computed without any difficulty.

\section{Summary}

We have extended the class of orthogonal polynomials for which the
existence of a commuting banded matrix for  the ``time-and-band limiting'' is
exhibited and then taken advantage of for numerical purposes. The search for
more instances of this important algebraic/numerical phenomenon remains a challenging
problem.

\section*{Acknowledgements}

The authors would like to thank Kenier Castillo and Lucas Tertuliano da Silva  for useful  comments on an  earlier version of the paper. 

\section*{Declarations}

\noindent\textbf{Use of artificial intelligence.} The assistance  of Claude (Anthropic) in this work is described in Sections~\ref{GenHermite} and~\ref{GenGegenbauer}. In particular, the argument leading up to and including Theorem~\ref{thmeig} was streamlined with the help of Claude. The authors have checked all the results and take full responsibility for the content of the manuscript.

\end{document}